\documentclass[%
 aip,
 amsmath,amssymb,
 reprint,%
]{revtex4-2}

\usepackage{graphicx}
\usepackage{dcolumn}
\usepackage{todonotes}

\usepackage{etoolbox}
\usepackage{epstopdf}
\usepackage[export]{adjustbox}

\usepackage{algorithm,algorithmic}
\usepackage{array} 
\usepackage{booktabs}
\usepackage{subfig}
\usepackage{amsmath}
\usepackage{amssymb}
\usepackage{blindtext}
\usepackage{xcolor}
\usepackage{soul}

\def\bs{\boldsymbol}

\newcolumntype{C}[1]{>{\centering}m{#1}}

\AtBeginDocument{%
  \providecommand\BibTeX{{%
    \normalfont B\kern-0.5em{\scshape i\kern-0.25em b}\kern-0.8em\TeX}}}

\makeatletter
\def\@email#1#2{%
 \endgroup
 \patchcmd{\titleblock@produce}
  {\frontmatter@RRAPformat}
  {\frontmatter@RRAPformat{\produce@RRAP{*#1\href{mailto:#2}{#2}}}\frontmatter@RRAPformat}
  {}{}
}%
\makeatother
\begin{document}

\preprint{AIP/123-QED}
\title[\resizebox{7in}{!}{Recurrent Neural Networks for Dynamical Systems: Applications to Ordinary Differential Equations, Collective Motion, and Hydrological Modeling}]{Recurrent Neural Networks for Dynamical Systems: Applications to Ordinary Differential Equations, Collective Motion, and Hydrological Modeling}

\author{Y. Park}

\author{K. Gajamannage}%
 \email{kelum.gajamannage@tamucc.edu}
\affiliation{ 
Department of Mathematics \& Statistics, Texas A\&M University -- Corpus Christi, Corpus Christi, TX 78412, USA
}%

\author{D.I. Jayathilake}%
\affiliation{ 
Department of Physical \& Environmental Sciences, Texas A\&M University -- Corpus Christi, Corpus Christi, TX 78412, USA
}%

\author{E. M. Bollt}
\affiliation{%
Department of Electrical and Computer Engineering, Clarkson University, Potsdam, NY 13699, USA 
}%

\date{\today}

\begin{abstract}
Classical methods of solving spatiotemporal dynamical systems include statistical approaches such as autoregressive integrated moving average, which assume linear and stationary relationships between systems' previous outputs. Development and implementation of linear methods are relatively simple, but they often do not capture non-linear relationships in the data. Thus, artificial neural networks (ANNs) are receiving attention from researchers in analyzing and forecasting dynamical systems. Recurrent neural networks (RNN), derived from feed-forward ANNs, use internal memory to process variable-length sequences of inputs. This allows RNNs to applicable for finding solutions for a vast variety of problems in spatiotemporal dynamical systems. Thus, in this paper, we utilize RNNs to treat some specific issues associated with dynamical systems. Specifically, we analyze the performance of RNNs applied to three tasks: reconstruction of correct Lorenz solutions for a system with a formulation error, reconstruction of corrupted collective motion trajectories, and forecasting of streamflow time series possessing spikes, representing three fields, namely, ordinary differential equations, collective motion, and hydrological modeling, respectively. We train and test RNNs uniquely in each task to demonstrate the broad applicability of RNNs in reconstruction and forecasting the dynamics of dynamical systems.
\end{abstract}

\maketitle
\begin{quotation}
Learning non-linear systems in nature requires researchers to solve complicated mathematical equations which often involve both high computational complexity and less computational precision \cite{hu2016new}. Artificial neural networks (ANNs) are a class of algorithms that model underlying relationships of a set of data through a process that mimics the way the human brain operates. Data-driven approaches such as ANNs are capable of affording such nonlinear systems with both high efficiency and good accuracy. ANN models are increasingly popular for two reasons. First, they guarantee that arbitrary continuous dynamical systems can be approximated by ANNs with a sufficient amount of sub-models which are called hidden layers \cite{hornik1989multilayer}. Second, ANNs themselves can be viewed as discretizations of continuous dynamical systems which make them suitable for studying dynamics \cite{wang1992discrete}. A recurrent neural network (RNN) is a class of ANNs where the connections between nodes form a directed graph along a temporal sequence where the output of a hidden node at the current time-step is sent to the corresponding hidden node for the next time-step. Among all the ANN architectures, RNNs are considered to be the most faithful candidates for modeling and forecasting temporal data sequences such as time series data. Thus, here we implement RNN to reconstruct and forecast data generated by several nonlinear dynamical systems.
\end{quotation}

\section{Introduction}\label{sec:level1}
Dynamical systems may exhibit dynamics that are highly sensitive to initial conditions, which is referred to as the butterfly effect. As a result of this sensitivity, which often manifests itself as an exponential growth of perturbations in the initial conditions, the behavior of dynamical systems is sophisticated \cite{jayawardena2013environmental}. From a historical point of view,  reconstruction or forecast of dynamical systems, especially those with spatiotemporal solutions, has been based on the techniques like linear parametric autoregressive models, moving-average models, autoregressive moving-average models, etc. All these models are linear and some of them possess the assumption that the data is being stationary, thus, these models are not able to cope with complex relationships in the data \cite{diaconescu2008use}.

\begin{figure}[htp]
\includegraphics[width = .9\linewidth]{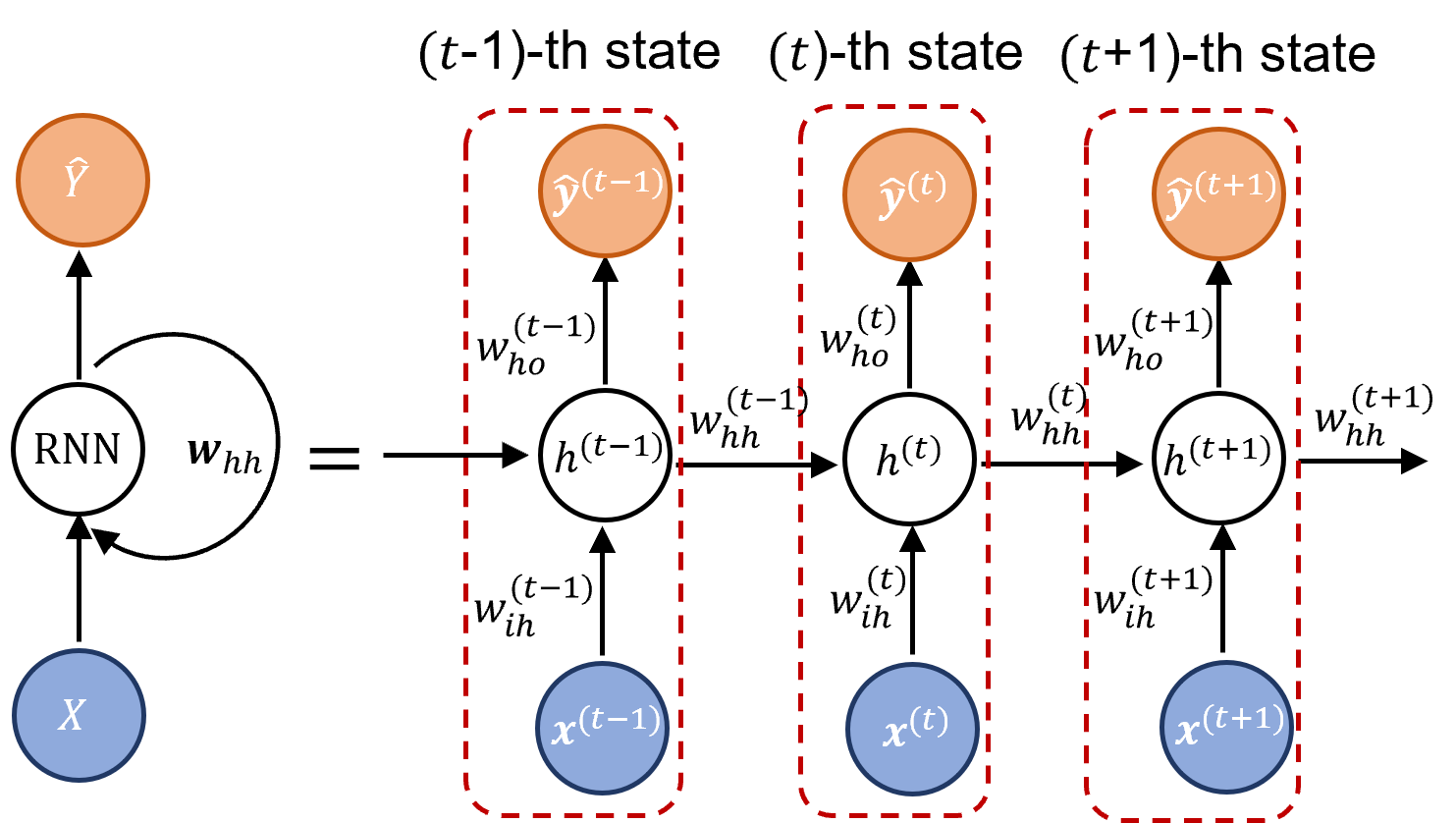}
\caption{Standard single-layer recurrent neural network (RNN), i.e., one-stacked RNN, architecture. The right-hand side schematic is the unrolled version, with respect to time, of the left-hand side RNN architecture where the directed loop represents recurrent feedback between different states of the same layer. Here, $X$, $\bs{w}_{hh}$, and $\hat{Y}$ represent input, the weights of the recurrent connections, and the output. Moreover, for the $t$-th state, while $\bs{h}^{(t)}$ represents the hidden state, $w^{(t)}_{ih}$, $w^{(t)}_{hh}$, and $w^{(t)}_{ho}$ represent weights between input and hidden layer, states of the hidden layer, and hidden layer and output, respectively. Note that, each state itself is a conventional artificial neural network; however, the RNN is formed by linking consecutive states within each layer with a directed edge. This recurrent feedback accords memory for an RNN  through its states (i.e., time-steps).}
\label{fig:rnn1}
\end{figure}

We utilize RNNs with unique ways of training and testing routings for modeling a general class of dynamical systems since they are considered to be the state-of-the-art machine learning methods for sequential data. Applications of RNNs include image classification, object recognition, object detection, speech recognition, language translation, voice synthesis, etc \cite{yuan2019adversarial}. An ANN usually involves a large number of parallelly operating processing units that are arranged in layers. The first layer receives the raw input information analogous to optic nerves in human visual processing. Each successive layer receives the output from both the layer preceding it and the same layer, rather than from the raw input, as analogous to the way that the neurons further from the optic nerve receive signals from those closer to it \cite{yuan2019adversarial}. The last layer produces the output of the system. RNN possesses an internal state or short-term memory due to the recurrent feedback connections, as seen by directed loops in Fig.~\ref{fig:rnn1}. This makes RNN is suitable for modeling sequential data, especially time series. Moreover, these connections allow well-trained RNNs to regenerate dynamics of any temporal nonlinear dynamical system up to a satisfactory accuracy. Thus, RNN models are widely used to analyze diverse time series problems including dynamical systems  \cite{sagheer2019unsupervised}. Most of the complex RNN architectures, such as Long Short-Term Memory (LSTM) \cite{hochreiter1997long} and Gated Recurrent Unit \cite{cho2014properties}, can be interpreted as a variation or as an extension of the basic RNN scheme\cite{bianchi2017overview}.

The objective of this study is to demonstrate how well RNNs learn spatiotemporal dynamics of multi-scale dynamical systems. Our focus is either reconstruct or forecast short-term or long-term trajectories of a dynamical system of interest with given initial conditions. We use three exemplary datasets that are generated from three dynamical systems, namely, the Lorenz attractor \cite{tabor1989chaos}, a generalized Vicsek model \cite{Gajamannage2019Reconstruction}, and a streamflow model, where those represent three diverse fields, ordinary differential equations, collective motion, and hydrological modeling, respectively. For given orbits that are generated from a Lorenz system with a formulation error, we train an RNN to eliminate the formulation error and reconstruct the correct system's responses. Then, we generate a collective motion dataset using the generalized version in \cite{Gajamannage2019Reconstruction}, of the classic Vicsek model in \cite{Vicsek1995}, and create two copies of that, one by imposing some additive Gaussian noise and the other by deleting some segments of the trajectories. Here, our goals in these experiments  are to eliminate the additive noise from the noisy trajectories. Finally, we forecast streamflow data that are naturally highly volatile and consist of frequent spiky fluctuations.

This paper is structured as follows: In Sec.~\ref{sec:rnn}, we present the architecture and the optimization routine of RNNs for analyzing and forecasting data generated from spatiotemporal dynamical systems.  In Sec.~\ref{sec:perfor}, we validate the performance of RNNs applied to three dynamical system datasets representing three fields, ordinary differential equations, collective motion, and hydrological modeling. The analysis conducted using RNNs addresses three unique issues and solutions for them. Sec.~\ref{sec:concl} states the discussion including key findings, limitations, future work, and conclusions. Table~\ref{tab:nom} represents notations used in this paper along with their descriptions.

\begin{table*}[htp]
\caption{Notations used in this paper and their descriptions}
\begin{tabular}{p{2cm}|p{13cm}}
Notation &  Description \\
\hline\\
$N_i$ & Number of nodes in the input layer of an RNN.\\
$N_h$ &  Number of nodes in the hidden layer of an RNN. \\
$N_o$ &  Number of nodes in the output layer of an RNN. \\ 
$N$ & Total number of orbits, trajectories, and agents.\\
$K$ & Total number of layers.\\
$t$ & Index for time-steps where $1\le t\le T.$\\
$k$ & Index for the RNN's layers where $1\le k\le K.$\\
$n$ & Index for orbits, trajectories, or agents where $1\le n\le N.$\\
$\sigma, \rho, \beta$ & Parameter in Lorenz.\\
$T$ & Number of total time-steps\\
$\delta$ & Step size of the generalized vicsek model\\
$\eta$ & Corruption level of orbits.\\
$\sigma$ & Standard deviation of Gaussian noise distribution.\\
$t$ & Index for time-steps where $1\le t\le T.$\\
$\epsilon^{(t)}_i$ & Noise parameter imposed to the $i$-th agent at the $t$-th time-step.\\
$\theta^{(t)}_i$ & Orientation of the $i$-th agent at the $t$-th time-step.\\
$r_d$ & Radius of interaction\\
$r^{(t)}$ & Radius of spiral example. \\
$k^{(t)}$ & Angle of rotation of the a spiral at time $t$ from the initial position. \\
$N^{(t)}_i$ & $r$-radius of neighborhood of the $i$-th agent at the $t$-th time-step.\\

\hline
$\bs{x}^{(t)}$ &  Input vector at time-step t.\\
$\bs{h}_k^{(t)}$ &  Hidden state of the $k$-th layer for the $t$-th time-step. \\
$\bs{y}^{(t)}$ &  Label vector at time-step t. \\
$\bs{b}_i$ &  Bias vector for the input layer.\\
$\bs{b}_h$ &  Bias vector for the hidden layer. \\
$\bs{x}^{(t)}_n$ &  Input vector of $n$-th orbit at time-step $t$ ($t$-th point of $n$-th orbit).\\
$\bs{y}^{(0)}_n$ &  Initial condition of $n$-th orbit.\\
$\bs{y}^{(t)}_n$ &  Label vector of the $n$-th orbit at time-step t. \\
$\hat{\bs{y}}^{(t)}_n$ &  Output vector of the $n$-th orbit at time $t$ corresponding to $\bs{x}^{(t)}_n$. \\
$\bs{u}^{(t)}_i$ & Average direction of motion of the agents in the neighborhood of $N^{(t)}_i$.\\
$\bs{c}^{(t)}$ & The $t$-th point on the centroid of the spiral. \\
$\bs{v}^{(t)}_i$ & Velocity of the $i$-th agent at the $t$-th time-step.\\

\hline
$X$ & Set of alll input data.\\
$\hat{Y}$ & Set of all output.\\
$Y$ & Set of all label.\\

\hline
$f$ & Linear transformation in RNN. \\
$g$ & Non-linear transformation in RNN. \\
$L(\cdot)$ & Loss function.\\

\hline
$W_{ih}$ & Weight matrix between the input layer $i$ and the hidden layer $h$.\\
$W_{h_k h_{k'}}$ & Weight matrix between the hidden layer $h_k$ and the hidden layer $h_{k'}$.\\
$W_{ho}$ & Weight matrix between the hidden layer $h$ and the output layer $o$.\\
$R^{(t)}_j$ & Rotation matrix of the $j$-th agent at the $t$-th time-step.\\
\end{tabular}\label{tab:nom}
\end {table*}

\section{Recurrent neural networks}\label{sec:rnn}
\begin{figure*}[hpt]
\includegraphics[width = .6\linewidth]{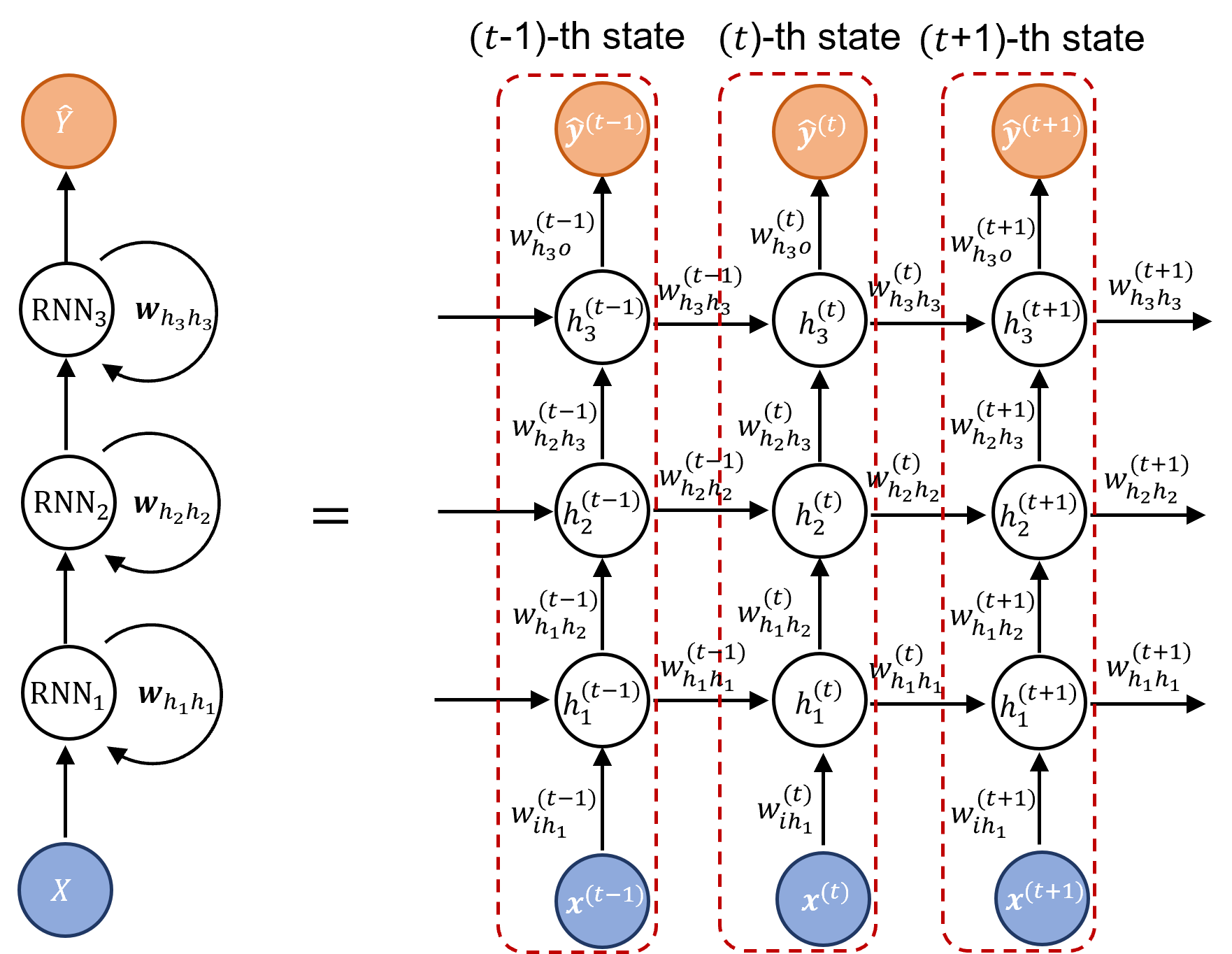}
\caption{Standard stacked recurrent neural network (RNN) architecture of having three layers, i.e., three-stacked RNN. The right-hand side schematic is the unrolled version, with respect to time, of the left-hand side RNN architecture where directed loops represent recurrent feedbacks between different states of the same layer. Here, $X$, $\bs{w}_{h_kh_k}$, and $Y$ represent input, the weights of the recurrent connections of the $k$-th layer, and the output. For the $t$-th state, where $1\le t\le T$ for some $T$, $\bs{x}^{(t)}$, $\bs{h}_k^{(t)}$, and $\hat{\bs{y}}^{(t)}$ represent the input vector, the hidden state of the $k$-th layer where $1\le k \le K$ for some $K$, and the output vector, respectively. Moreover, $w^{(t)}_{ih_1}$, $w^{(t)}_{h_kh_{k'}}$, and $w^{(t)}_{h_3o}$ are the weights between input and hidden layer 1, hidden layer $h_k$ and hidden layer $h_{k'}$, and hidden layer 3 and output, respectively. For the $k$-th layer, $w^{(t)}_{h_k h_k}$ represents the weight between two consecutive states.}
\label{fig:rnn2}
\end{figure*}

RNNs are a major class of machine learning tools that are particularly well suited to work on sequential data. Connections between nodes form a directed graph along a temporal sequence, that can be used to learn nonlinear dependencies of a system over time \cite{csaji2001approximation}. In its most general form, an RNN can be seen as a weighted  directed graph that contains three different kinds of nodes, namely, the input nodes, hidden nodes, and output nodes \cite{zhang2016architectural}. As seen in Fig.~\ref{fig:rnn1}, where we represent the structure of a simple RNN with one layer, input nodes do not have incoming connections and output nodes do not have outgoing connections, but hidden nodes have both. By RNN's design, two different nodes which are either at the same time level or at different time levels can be connected by an edge, as seen in Fig.~\ref{fig:rnn2}. 

An RNN has a memory that stores \emph{knowledge} about the data sequence that it has already seen, however, its memory of internal states through iterations is short-term \cite{bengio1994learning}. By applying activation functions (also called, transfer functions) to previous states and inputs, RNNs compute new states. Typical activation functions in RNNs include sigmoid, hyperbolic tangent, and ReLU; however, sigmoid and hyperbolic tangent encounter vanishing of the gradient in the optimization step \cite{chollet2017}. Properly customized RNNs with the best activation function and the best number of hidden nodes can capture the dynamic of any nonlinear system up to a high precision \cite{siegelmann1991turing}. The model capacity and flexibility for learning non-linear systems can be increased by a multi-layer recurrent structure \cite{schmidhuber1992learning}. By this multi-layer structure, the hidden output of one recurrent layer can be propagated through time and is used as the input data to the next recurrent layer. In our study, we used three-stacked RNNs (or recurrent layers) as seen in Fig.~\ref{fig:rnn2}. 

\subsection{Architecture of RNNs}
For simplicity, we formulate an RNN with only one hidden layer, i.e., one-stacked RNN, as shown in Fig.~\ref{fig:rnn1}. We denote the number of nodes in the input layer, hidden layer, and output layer as $N_i$, $N_h$, and $N_o$. For the $t$-th state, the weight matrix from input to hidden layer is denoted as $w^{(t)}_{ih} \in \mathbb{R}^{N_i \times N_h}$, the weight matrix between hidden layer as $w^{(t)}_{hh} \in \mathbb{R}^{N_h \times N_h}$, and the weight matrix from the hidden layer to the output layer as $w^{(t)}_{ho} \in \mathbb{R}^{N_h \times N_o}$. Moreover, for the $t$-th state, we denote the bias vector for the input layer as $\bs{b}^{(t)}_i \in \mathbb{R}^{N_i}$ and the bias vector for the hidden layer as $\bs{b}^{(t)}_h \in \mathbb{R}^{N_h}$. Here, we develop the formulation for the $t$-th state, where $1\le t\le T$ for some $T$, i.e. for the $t$-th ANN, in the RNN. The input $\bs{x}^{(t)}$ of the $t$-th state is multiplied with $w^{(t)}_{ih}$, and summed that with  $\bs{b}^{(t)}_i$, i.e. $w^{(t)}_{ih}\bs{x}^{(t)} +\bs{b}^{(t)}_i$, to get the flow of information from the input to the hidden layer. Similarly, the information flow by the recurrent feedback from the hidden node of the ($t-1$)-th state to that of the $t$-th state is $w^{(t-1)}_{hh}\bs{h}^{(t-1)} +\bs{b}^{(t-1)}_h$. The hidden state $\bs{h}^{(t)}$ of the $t$-th state is the sum of the above two information flows passing through a nonlinear transformation, say $f$, such that,
 \begin{equation}\label{eq:rnn1}
\bs{h}^{(t)} = f\left(w^{(t)}_{ih} \bs{x}^{(t)} + \bs{b}^{(t)}_i + w^{(t-1)}_{hh} \bs{h}^{(t-1)} + \bs{b}^{(t-1)}_h\right),
\end{equation}
where $f$ is the activation function of the neurons which is usually sigmoid or hyperbolic tangent. The information flow from the hidden layer of the $t$-th state to its output, i.e. $w^{(t)}_{ho}\bs{h}^{(t)}+\bs{b}_o$, is composed with another transformation, say $g$, to get the output, say $\hat{\bs{y}}^{(t)}$, of the RNN such that,
 \begin{equation}\label{eq:rnn2}
\hat{\bs{y}}^{(t)} = g\left(w^{(t)}_{ho} \bs{h}^{(t)} + \bs{b}^{(t)}_o\right),
\end{equation}
where $g$ is linear in general. For all $t$'s, $\bs{h}^{(t)}$ is usually initialized with a vector of zeros.  The explicit correspondence of Eqns.~\eqref{eq:rnn1} and \eqref{eq:rnn2} with the RNN's architecture, helps us adopt spatiotemporal notion of RNNs using so-called \emph{finite unfolding in time} with shared weight matrices \cite{schafer2007recurrent}.

By learning the formulation pattern in Eqns.~\eqref{eq:rnn1} and \eqref{eq:rnn2}, we derive the formulation for a stacked RNN with $K$ hidden layers. Note that, Fig.~\ref{fig:rnn2} provides the architecture for such an RNN with $K=3$, i.e., three-stacked RNN. For the $t$-th time-step, let $\bs{x}^{(t)}$, $\bs{h}_k^{(t)}$, and $\bs{y}^{(t)}$ represent the input vector, the hidden state of the $k$-th layer, and the output vector, respectively. We denote the number of nodes in the input layer, $k$-th hidden layer, and output layer as $N_i$, $N_{h_k}$, and $N_o$. Moreover, $w^{(t)}_{ih_1} \in \mathbb{R}^{N_i \times N_{h_1}}$, $w^{(t)}_{h_kh_{k'}}\in \mathbb{R}^{N_{h_k} \times N_{h_k'}}$, and $w^{(t)}_{h_3o}\in \mathbb{R}^{N_{h_3} \times N_o}$ are the weights between the input and hidden layer 1, hidden layer $h_k$ and hidden layer $h_{k'}$, and hidden layer 3 and the output, respectively. Similar to Eqn.~\eqref{eq:rnn1}, the $t$-th state of the first hidden layer of the RNN is sent through a nonlinear activation function $f$ along with both the input $\bs{x}^{(t)}$ and the $(t-1)$-th state of the first hidden layer, $\bs{h}^{(t-1)}_1$, as
\begin{equation} \label{eq:stacked1}
\bs{h}^{(t)}_1 = f\left(w^{(t)}_{ih_1} \bs{x}^{(t)} + \bs{b}^{(t)}_i + w^{(t-1)}_{h_1h_1} \bs{h}_1^{(t-1)} + \bs{b}^{(t-1)}_{h_1}\right),
\end{equation}
where $\bs{b}^{(t)}_i$ is the bias vector for the input at the $t$-th state and $\bs{b}^{(t-1)}_{h_1}$ is the bias vector for the first hidden layer at the $t$-th state. The hidden state, $\bs{h}^{(t)}_k$ for $k=2,3,\dots, K$ of the $t$-th state is based on both its $t$-th state at the previous hidden layer, i.e., $\bs{h}^{(t)}_{k-1}$, and the $(t-1)$-th state of the same hidden layer, i.e., $\bs{h}^{(t-1)}_k$, as
\begin{equation} \label{eq:stacked2}
\bs{h}^{(t)}_k = f\left(w^{(t)}_{h_{k-1}h_{k}} \bs{h}^{(t)}_{k-1} + \bs{b}^{(t)}_{h_{k-1}h_{k}} + w^{(t-1)}_{h_kh_k} \bs{h}^{(t-1)}_k + \bs{b}^{(t-1)}_{h_kh_k}\right),
\end{equation}
where $\bs{b}^{(t)}_{h_{k-1}h_{k}}$ is the bias vector for the $(k-1)$-th hidden later at the $t$-th state and $\bs{b}^{(t-1)}_{h_kh_k}$ is the bias vector for the $k$-th hidden layer at the $(t-1)$-th state. The output of the $t$-th state, say $\hat{\bs{y}}^{(t)}$, is only based on its $K$-th hidden layer such that,
\begin{equation} \label{eq:stacked3}
\hat{\bs{y}}^{(t)} = g\left(w^{(t)}_{h_Ko} \bs{h}^{(t)}_{K} + \bs{b}^{(t)}_{h_Ko}\right),
\end{equation}
where $\bs{b}^{(t)}_{h_Ko}$ is the bias vector for the $K$-th hidden layer at the $t$-th state. The stacked RNN's initial hidden states $\bs{h}^{(0)}_k,\dots,\bs{h}^{(T)}_k$ should be initialized for each layer $k$. All the weight matrices and bias vectors are optimized using the gradient descent method according to a straightforward backpropagation through time (BTT) procedure \cite{werbos1990backpropagation}.

\subsection{Optimization of RNNs}\label{sec:opt}
Training an RNN is mathematically implemented as a minimization of a relevant reconstruction error, widely called as the loss, function with respect to weights and bias vectors of Eqns.~\eqref{eq:stacked1}, \eqref{eq:stacked2}, and \eqref{eq:stacked3}. This optimization is carried out in four steps: first, forward propagation of input data through the neural network to get the output; second, calculate the loss between forecasted output and the expected output; third, calculate the derivatives of the loss function with respect to the ANN's weights and bias vectors using BTT; and fourth, adjusting the weights and bias vectors by gradient descent method \cite{gruslys2016memory}.

BTT unrolls backward all the dependencies of the output on the weights of the system \cite{manneschi2020alternative}, which is represented from left side to right side in Figs.~\ref{fig:rnn1} and \ref{fig:rnn2}. We train the RNN by one data instance, say instant $n$, at a time; thus, we formulate the loss function of RNN's as a normalization over all the data instances, say $N$. We denote the $n$-th input data instance as $\bs{x}_n =[\bs{x}^{(1)}_n, \dots, \bs{x}^{(T)}_n] \in X$ where $X$ represents the set of all the input data, and the corresponding output as $\hat{\bs{y}}_n =[\hat{\bs{y}}^{(1)}_n, \dots, \hat{\bs{y}}^{(T)}_n]\in \hat{Y}$ where $\hat{Y}$ represents the set of all the outputs. We denote the label of the $n$-th input data as $\bs{y}_n = [\bs{y}^{(1)}_n, \dots, \bs{y}^{(T)}_n]\in Y$ where $Y$ represents the set of all the labels for the input data. The loss function, denoted by $L$, is defined as the normalized squared difference between labels and RNN's output,
\begin{equation}\label{eq:mse}
L \left(\hat{Y}, Y\right) = \frac{1}{|Y|} \sum_{ \bs{y}_n \in Y } ( \hat{\bs{y}}_n - \bs{y}_n )^2  
\end{equation}
where the cardinality $|Y|$ is defined as the number of data instances in the set $Y$. We use BTT to compute the derivatives of Eqn.~\eqref{eq:mse} with respect to the weights and bias vectors. We update the weights using the gradient descent based method, called Adaptive Moment Estimation (ADAM) \cite{kingma2014adam}. ADAM is an iterative optimization algorithm that is widely used in modern machine learning algorithms to minimize loss functions where it employs the averages of both the first moment gradients and the second moment of the gradients for computations. ADAM generally converges faster than standard gradient descent methods and saves memory by not accumulating the intermediate weights.

\section{Performance analysis}\label{sec:perfor}
In this section, we work on three diverse problems: 1) a Lorenz system having a formulation error, 2) denoising of noisy trajectories of particle swarms in collective motion, and 3) forecasting of groundwater streamflow of hydrologic catchments. These three problems, in that order, represent three fields, ordinary differential equations, collective motion, and hydrological modeling.

Here, we use the RNN model available in the machine learning library \emph{Pytorch} \cite{NEURIPS2019-9015}, which requires the parameter inputs sequential length, input size, learning rate, output size, number of RNNs, and number of hidden nodes. The sequential length is the length of the input data sequence, which is $T$ in Fig.~\ref{fig:rnn2}. The input size is the number of features or dimensions in the dataset, which is the length of the vector $\bs{x}^{(t)}$ in Fig.~\ref{fig:rnn2}. The learning rate is associated with the ADAM optimization routing where a large learning rate allows the RNN to learn faster at a cost of arriving into a sub-optimal final set of weights. A smaller learning rate may allow the RNN to learn further optimal or even globally optimal set of weights but may take significantly longer time to train. The number of RNNs is the number of basic RNN units that are stacked together which is equal to the number of hidden layers, where it is $K=3$ in Fig.~\ref{fig:rnn2}. The number of hidden nodes is the user input number of nodes at each hidden layer of the RNN, which is $N_{h_k}$ for $k$-th hidden layer of Fig.~\ref{fig:rnn2}. Finding both the optimal number of RNNs and the optimal number of hidden nodes could prevent possible over-fitting or under-fitting scenarios of the RNN's training process. 

\subsection{Lorenz system with a formulation error}
The Lorenz system is a chaotic dynamical system \cite{devaney2018first}, that is formulated as a set of ordinary differential equations which is known for having chaotic solutions for certain parameter values and initial conditions \cite{lorenz1963deterministic}. Here, we assume that we are given a Lorenz system with a formulation error so that our task is to make an ANN to generate the true solutions related to the input erroneous solutions generated by the system with the formulation error. The standard Lorenz system of the temporal variable $t$ and three spacial variables $y_1$, $y_2$, and $y_3$ is given by 
\begin{equation}\label{eq:lorenz1}
\begin{aligned}
\frac{dy_1}{dt} = \sigma(y_2-y_1), \\
\frac{dy_2}{dt} = y_1(\rho - x_3) - y_2, \\
\frac{dy_3}{dt} = y_1y_2 - \beta y_3, 
\end{aligned}
\end{equation}
where $\sigma$, $\rho$, and $\beta \in \mathbb{R}$ are parameters. To mimic a formulation error, we eliminate the term $-y_2$ from Eqn.~\eqref{eq:lorenz1} as
\begin{equation}
\label{eq:lorenz2}
\begin{aligned}
\frac{dx_1}{dt} = \sigma(x_2-x_1), \\
\frac{dx_2}{dt} = x_1(\rho - x_3), \\
\frac{dx_3}{dt} = x_1x_2 - \beta x_3. 
\end{aligned}
\end{equation}
Note that, we replaced the variable $y$ in Eqn.~\eqref{eq:lorenz1} with variable $x$ to distinguish data from two systems. We set $\sigma=3$, $\rho=28$, and $\beta=8/3$ in both Eqn.~\eqref{eq:lorenz1} and Eqn.~\eqref{eq:lorenz2} since the solutions with these parameters are well-known for showing chaotic behaviors. For different initial conditions, we solve the system in Eqn.~\eqref{eq:lorenz1} and treat them as training labels for the RNN. Now, we generate the erroneous solutions for Eqn.~\eqref{eq:lorenz2} based on Eqn.~\eqref{eq:lorenz1} to use as training data as we will explain in the sequel. We use the trained RNN to generate the true solutions for the input erroneous solutions.

\begin{figure}[htp]
\centering
\includegraphics[width = 2.5in]{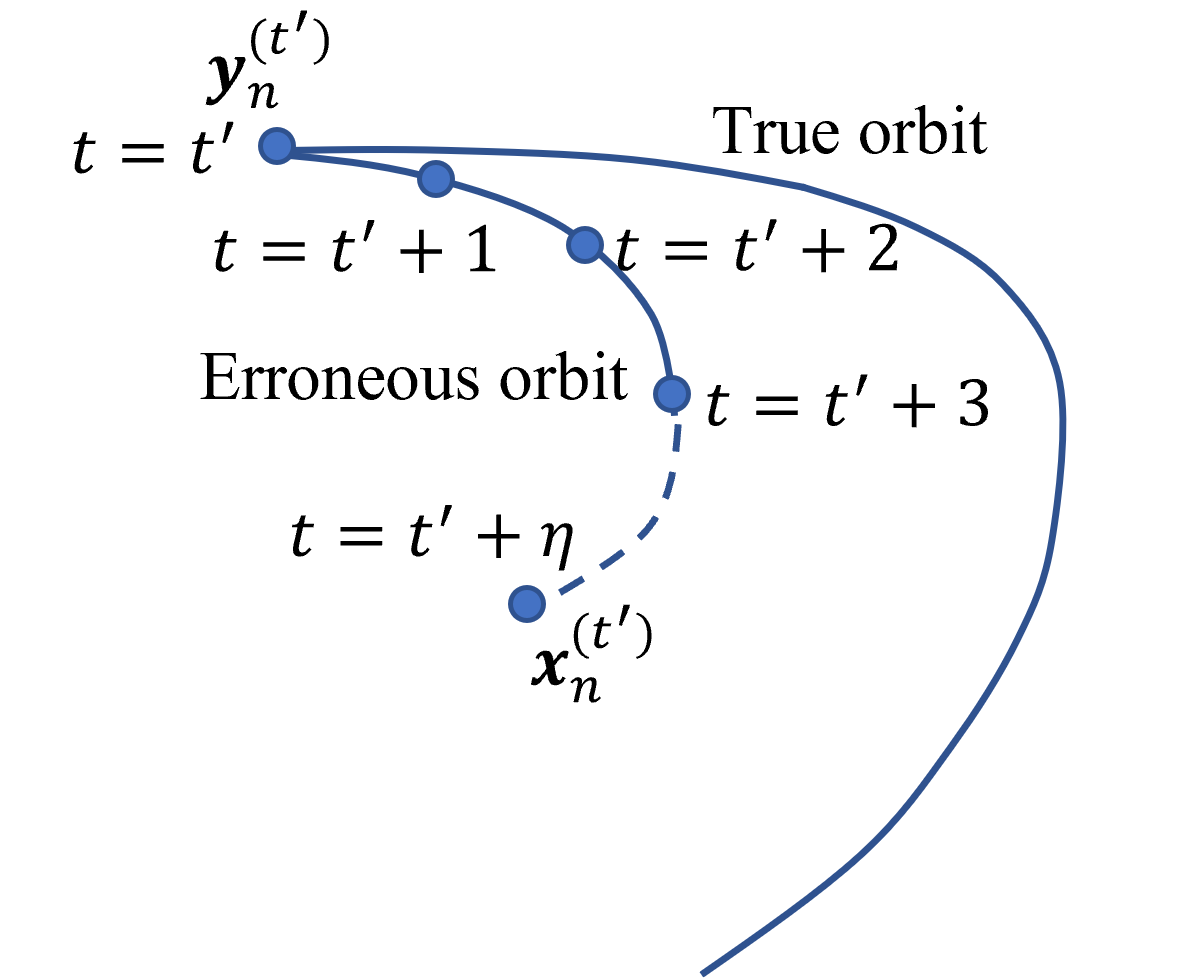}
\caption{For some $t^\prime$, where $1\le t^\prime\le T$, we treat $\bs{y}_n^{(t^\prime)}$ in the $n$-th true orbit as the initial condition for the erroneous Lorenz system and solve it for some time-steps, denoted as $\eta$, with the same parameter values used for solving the true Lorenz system. We consider the last point $\eta$, denoted by $\bs{x}_n^{(t^\prime)}$, of the solution set as the erroneous solution relevant to the true point $\bs{y}_n^{(t^\prime)}$. We carry out this process for each $t^\prime$ to generate the corrupted orbit $\left[\bs{x}_n^{(0)}\right.$, $\dots$, $\bs{x}_n^{(t)}$, $\dots$, $\left.\bs{x}_n^{(T)}\right]$. Here, $\eta$ can be considered as the corruption level of the orbit where a big $\eta$ contributes a high corruption.}
\label{fig:eta}
\end{figure}

We generate 500 random initial conditions, $\bs{y}_n^{(0)}=$ $\left(y_{n,1}^{(0)}, y_{n,2}^{(0)},y_{n,3}^{(0)}\right)$ where 1 $\le$ $n$ $\le$ 500, from the uniform distribution $\mathbb{U}[-15,15]$. We solve Eqn.~\eqref{eq:lorenz1} and generate 500 orbits, denoted as $N$ with the step-size, denoted as $\delta$, is 0.01, and the total time-steps, denoted as $T$, is 5000 (i.e., $t\in[0,50]$ of Eqn.~\eqref{eq:lorenz1}). Here, the $n$-th orbit is given as $\bs{y}_n=\left[\bs{y}_n^{(0)}\right.$, $\dots$, $\bs{y}_n^{(t)}$, $\dots$, $\left.\bs{y}_n^{(5000)}\right]$, where $\bs{y}_n^{(t)}= {\left(y_{n,1}^{(t)}, y_{n,2}^{(t)},y_{n,3}^{(t)}\right)}^T$. Now, we explain the process of generating the erroneous orbit $\bs{x}_n=\left[\bs{x}_n^{(0)}\right.$, $\dots$, $\bs{x}_n^{(t)}$, $\dots$, $\left.\bs{x}_n^{(5000)}\right]$ from the above true orbit. For some $t$, we treat $\bs{y}_n^{(t)}$ in the $n$-th orbit as an initial condition for the system in Eqn.~\eqref{eq:lorenz2} and solve it for some time-steps, denoted as $\eta$, where $\eta<T$, with the same values for the parameters $\sigma$, $\rho$, $\beta$, and $\delta$ as above. The parameter $\eta$ governs the corruption level of the orbits such that a big value of $\eta$ contributes a big level of corruption to the orbits. We treat the last point of the solution set, that is the $\eta$-th point, as the erroneous solution relevant to the correct point $\bs{y}_n^{(t)}$, that we denote by $\bs{x}_n^{(t)}$, see Fig.~\ref{fig:eta}. We carry out this process for each $n$ and generate the entire erroneous orbit $\left[\bs{x}_n^{(0)}\right.$, $\dots$, $\bs{x}_n^{(t)}$, $\dots$, $\left.\bs{x}_n^{(5000)}\right]$. Similarly, we compute the erroneous orbits corresponding to each of the 500 true orbits that we generated by changing initial conditions. The first 400 orbits of the system in Eqn.~\eqref{eq:lorenz2} are chosen as the training set, and the corresponding 400 orbits that are derived from the system in Eqn.~\eqref{eq:lorenz1} are chosen as their corresponding labels.

First, we set $\eta=1$ and generate 500 true and erroneous orbits each with the 5000 time-steps as explained above, where the first 400 orbits are used for training while the last 100 orbits are used for testing. In order to avoid both the over-fitting and the semi-convergence of RNN, the best number of training epochs was chosen by running the RNN twice as follows: first, we run the RNN with a fixed, but big, number of epochs to find the epoch number that gives the minimum loss; second, we re-run the RNN for that many epochs to complete the training. For the Lorenz system, based on the corruption level $\eta$, the optimum number of epochs ranges between 179 and 930. We utilize the activation function ReLU as it improves the reconstruction error about 10 times than that of the hyperbolic tangent activation function. In order to find the optimal parameters, which are learning rate, number of RNN, number of hidden node; first, we discretize the parameter space; then, we train and test the RNN with those discretized parameter combinations; finally, we select the best values for the parameters that give the least reconstruction error. In all the experiments, a data batch is equal to the entire dataset. We set the other parameters of RNN as follows: sequential length to 5000 since each orbit has 5000 time-steps; input size to three since each orbit is three-dimensional; output size to three since the output is also three-dimensional; the learning rate of the Lorenz orbits to 0.01.

\begin{figure*}[htp]
\includegraphics[width = 1\linewidth]{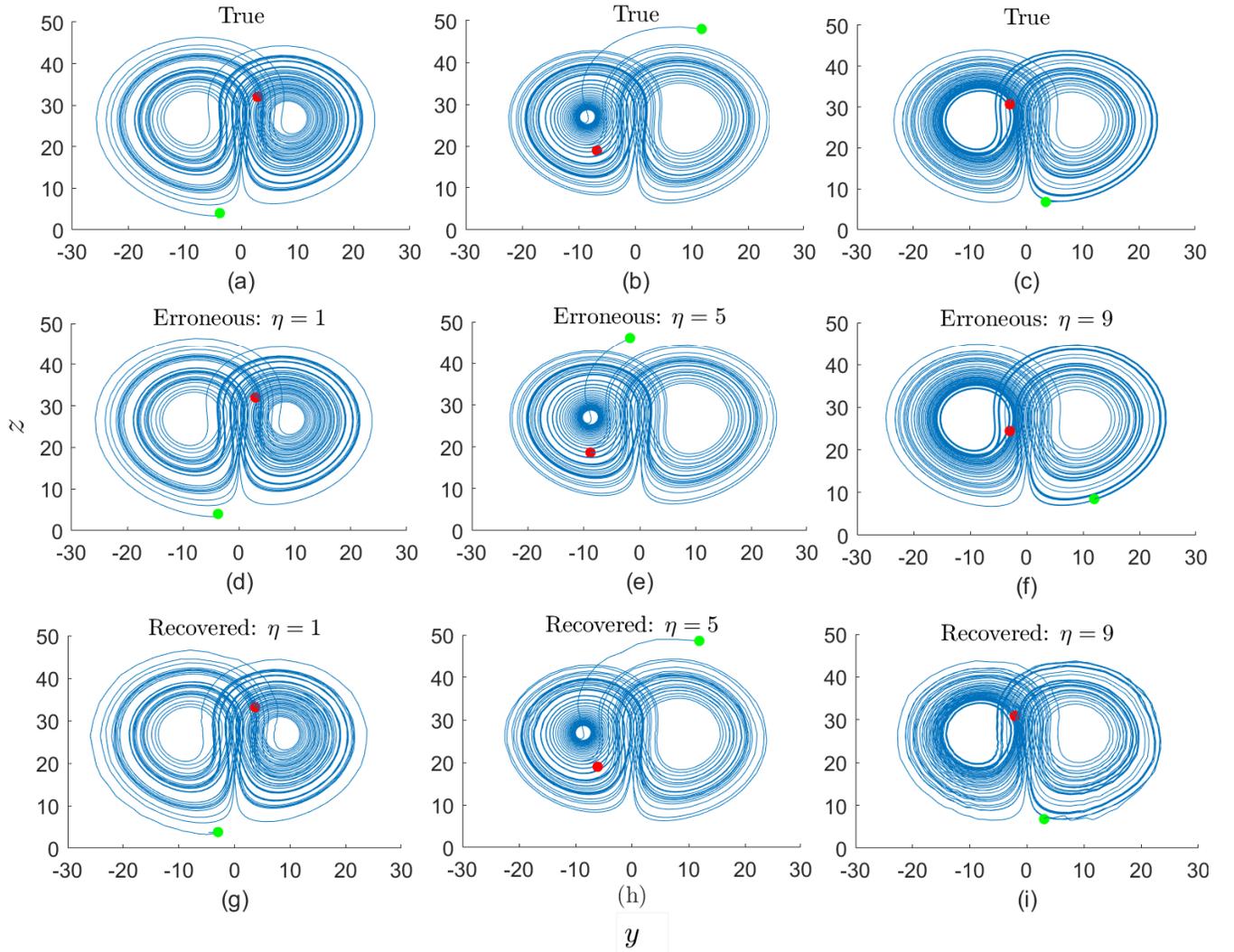}  
\caption{Recovery of the true solution of an erroneous Lorenz system using RNNs. Here, all the figures are two-dimensional ($y$-$z$ plan) projections of the three-dimensional orbits. 120 true orbits are generated using a Lorenz system and a Lorenz system with a formulation error is used to generate corrupted orbits corresponding the true orbits with some known corruption level $\eta$. Here, while (a-c) show three arbitrarily chosen true orbits, (d-f) show the erroneous orbits corresponding to them, respectively, where the green disks represent the starting point and red disks represent the ends. For each $\eta$, we set the parameters of the RNN as follows: 5000 for sequential length, 3 for both input size and output size, 0.01 for learning rate, three-stacked RNN with 128 number of hidden node each, and ReLu for the activation function. Moreover, the best number of epochs vary between 179 and 930 based on $\eta$. First, we train an RNN with arbitrary 100 orbits for each $\eta$, and then we pass the erroneous orbits through the RNNs and recovery them as show in (g)-(i). The start and the end points of the recovered orbits are close to those of the true orbits whereas they are far apart between erroneous and true orbits. Root mean square error of the reconstructions are 0.072, 0.15, 0.27 for the corruption levels $\eta = 1, 5$, and 9, respectively.}
\label{fig:orbit}
\end{figure*}

We set the RNN to three hidden recurrent layers with $256$ nodes each.  For each $1\le n\le 400$, we input each $\bs{x}_n$ into the RNN one at a time and generate the corresponding output of the RNN, denoted as $\hat{y}_n$. We compute the reconstruction error $L$ between the outputs and the labels using Eqn.~\eqref{eq:mse}. We use BTT to compute the derivatives of Eqn.~\eqref{eq:mse} with respect to the weights and the bias vectors, and update the weights and the bias vectors using ADAM as explained in Sec.~\ref{sec:opt}. The remaining 100 test orbits, $\bs{x}_n$ where $401\le n\le 500$, are fed into this trained RNN one orbit at a time and obtain the recovered solutions, $\hat{\bs{y}}_n^{(t)}; 401\le n\le 500$. We repeat the same experiment nine more times with $\eta=2, 3 \dots, 9$. Fig.~\ref{fig:orbit} shows two-dimensional projections of the three-dimensional true, erroneous, and recovered orbits for the cases $\eta$ = 1, 5, and 9.

\begin{figure}[htp]
\includegraphics[width = 3in]{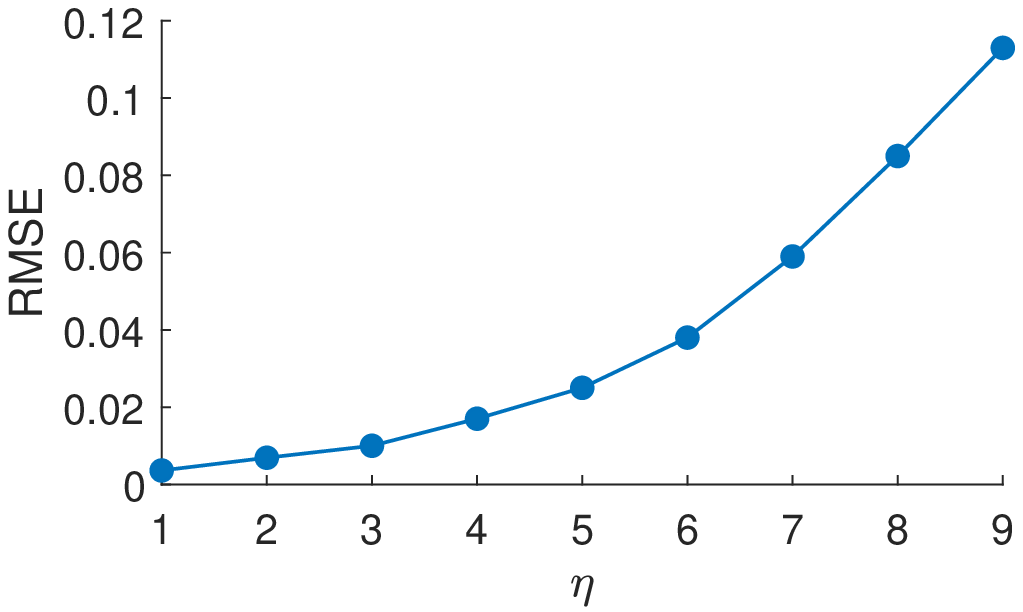}
\caption{RMSE of the reconstruction of corrupted Lorenz orbits versus different levels of the corruptions ($\eta$). 120 true orbits are constructed using a conventional Lorenz system by changing initial conditions. Then, these true orbits are used as the initial conditions to generate solutions for the Lorenz system with a formulation error where this generates 120 corrupted orbits. Each of these corrupted orbits has undergone $\eta$-level of corruption where $\eta=1, 2, \dots, 9$ to create 9 separate datasets each with 120 true orbits and 120 corrupted orbits. Then, for each $\eta$, an RNN having three hidden layers is trained with the first 100 orbits and then test on the remaining 20 orbits. The root mean squared errors (RMSE) between the true and the reconstructed trajectories justify that the reconstruction error increases rapidly when the corruption level increases.}
\label{fig:kLoss}
\end{figure}

We observe that due to the formulation error of the Lorenz system the erroneous orbits show some shifting of the initial and the end points from their true locations. This shift does not only occur at the two extremes but all the intermediate points are also shifted along the orbit in the same direction. We observed that the recovered orbits have their extreme points around their true locations. Moreover, the recovery quality gets worse as $\eta$ increases. To analyze the recovery performance with increasing corruption levels, we compute the normalized root mean square error, denoted by RMSE, between the recovered and the true orbits as
\begin{equation} \label{eq:rmse}
RMSE =\sqrt{\frac{1}{|Y|} \sum_{\forall \bs{y}_n\in Y}\|\bs{y}_n-\hat{\bs{y}}_n\|_F^2},
\end{equation} 
where $|Y|$ denotes the carnality, i.e, the number of observations in the test set, and  $F$ denotes the Frobenius norm. Fig.~\ref{fig:kLoss} shows that the RMSE of the recovery with respect to the corruption level where we see that the reconstruction error increases rapidly with respect to the increase of the corruption level.

\subsection{Noisy trajectories of collective motion}
Analysis of the trajectories of collectively moving agents such as fish \cite{gajamannage2015dimensionality} and birds \cite{gajamannage2015identifying}, is a highly active field in computer vision. Tracked trajectories of such agents are often contaminated with noise due to the causes like insufficient camera precision and lack of accuracy of the tracking method \cite{gajamannage2021reconstruction}. Here, we use RNNs to reconstruct the true collective motion trajectories from observed noisy trajectories. We simulate the collective motion trajectories using a modified version of a classic self-propelled particle model, named the Vicsec model \cite{Vicsek1995}.

The collective motion is defined as a spontaneous emergence of the ordered movement in a system consisting of many self-propelled agents. The Vicsek model given in Ref.~\onlinecite{Vicsek1995} is one of the widely used models that describe such behavior. This model performs simulations on a square-shaped region with periodic boundary conditions. We generate a synthetic collective motion dataset by using the generalized version presented in Ref.~\onlinecite{Gajamannage2019Reconstruction} of the classic Vicsek model where we incorporated a rotational matrix, denoted as $R^{(t)}_n$. The agent-wise temporal rotation matrix $R^{(t)}_n$ imposed to the $n$-th agent at the $t$-th time-step allows us to simulate interesting collective motion scenarios while ensuring the intra-group interactions between the agents. Based on the average direction of motion of all the particles in a neighborhood, denoted by $N^{(t)}_{n}$, within a radius $r_d$ of the $n$-th particle, the generalized Vicsek model updates the orientation of the $n$-th particle at the $t$-th time-step, denoted by $\theta_n^{(t)}\in[-\pi, \pi)$.
Therefore, the orientation of the $n$-th agent at the $t$-th time-step, defined as $\theta_n^{(t+1)}$, is computed as
\begin{equation}\label{eqn:modVicsek1}
\theta_n^{(t+1)}=\arg(\bs{u}^{(t)}_n)+\epsilon_n^{(t)},
\end{equation}
where, $\epsilon_n^{(t)}$ is a noise parameter imposed to the orientation of the $n$-th agent at the $t$-th time-step, and 
\begin{equation}\label{eqn:modVicsek}
\bs{u}^{(t)}_n=\frac{1}{|N^{(t)}_{n}|} \sum_{j\in N^{(t)}_{n}} R^{(t)}_j \begin{bmatrix}\cos(\theta_j^{(t)}) \\ \sin(\theta_j^{(t)})\end{bmatrix}.
\end{equation}
The position vector of the $n$-th agent at the $t$-th time-step is given as
\begin{equation}\label{eqn:modVicsek2}
\bs{x}_n^{(t+1)}=\bs{x}_n^{(t)}+v_n^{(t)}R_n^{(t)} \begin{bmatrix}\cos(\theta_n^{(t)}) \\ \sin(\theta_n^{(t)})\end{bmatrix} \delta, \\
\end{equation}
where, $v_n^{(t)}$ denotes the speed of the $n$-th agent at the $t$-th time-step and $\delta$ is the step size. 
Eqns.~\eqref{eqn:modVicsek1}, \eqref{eqn:modVicsek}, and \eqref{eqn:modVicsek2} with $R_n^{(t)}=[1,0;0,1]$ for all $n$'s and $t$'s implies the classic Vicsek model.

We use this rotational matrix to formulate a spiral collective motion scenario since spiral collective motion is often observed in nature ranging from biology, such as bacteria colonies \cite{Koizumi2020}, to astronomy, such as spiral galaxies \cite{Parnovsky2010}. Here, first, we formulate an anticlockwise Archimedean spiral that rotates an angle of $3\pi$ and then compute the agent-wise temporal rotational matrix based on that spiral. We assume that the two-dimensional coordinates $\left[\bs{c}^{(1)}; \dots; \bs{c}^{(t)}; \dots; \bs{c}^{(T)}\right]$ where $\bs{c}^{(t)} \in\mathbb{R}^2$ represent this spiral such that
\begin{equation}
\begin{split}
\bs{c}^{(t)}=r^{(t)} \begin{pmatrix} \cos\left(\kappa^{(t)}\right) \\ \sin\left(\kappa^{(t)}\right) \end{pmatrix} \ ; \ \text{where} \\
r^{(t)} = 1+\frac{3}{T-1}(t-1) \ \text{and} \ \kappa^{(t)}=\frac{3\pi}{T}(t-1),
\end{split}
\end{equation}
for $t=1, \dots, T$. Here, $r^{(t)}$ is the variable radius of the spiral where it changes from 1 to 4 and $\kappa^{(t)}$ is the angle of rotation with respect to the origin that varies from 0 to $3\pi$. The rotation angle with respect to the $(t-1)$-th coordinates, $\bs{c}^{(t-1)}=\left[c^{(t-1)}_1, c^{(t-1)}_2\right]^{T_r}$, of this spiral, denoted by $\gamma^{(t)}$, is
\begin{equation}
\gamma^{(t)}=\tan^{-1}\left(\frac{c_2^{(t)}-c_{2}^{(t-1)}}{c_1^{(t)}-c_{1}^{(t-1)}}\right) \ ; \ t=2,\dots,T.
\end{equation} 
Thus, the two-dimensional rotational matrix for the $i$-th agent at the $t$-th time-step is given by
\begin{equation}\label{eqn:rotmat2}
R^{(t)}_{n}=  
\begin{pmatrix} \cos\left(-\gamma^{(t)}\right) & -\sin\left(-\gamma^{(t)}\right) \\ \sin\left(-\gamma^{(t)}\right) & \cos\left(-\gamma^{(t)}\right) \end{pmatrix}
\end{equation}
Eqns.~\eqref{eqn:modVicsek1}, \eqref{eqn:modVicsek}, and \eqref{eqn:modVicsek2} with the rotational matrix in Eqn.~\eqref{eqn:rotmat2} provide the complete formulation of the system generating the spiral collective motion dataset.

\begin{figure*}
\centering
\includegraphics[width= .95\textwidth]{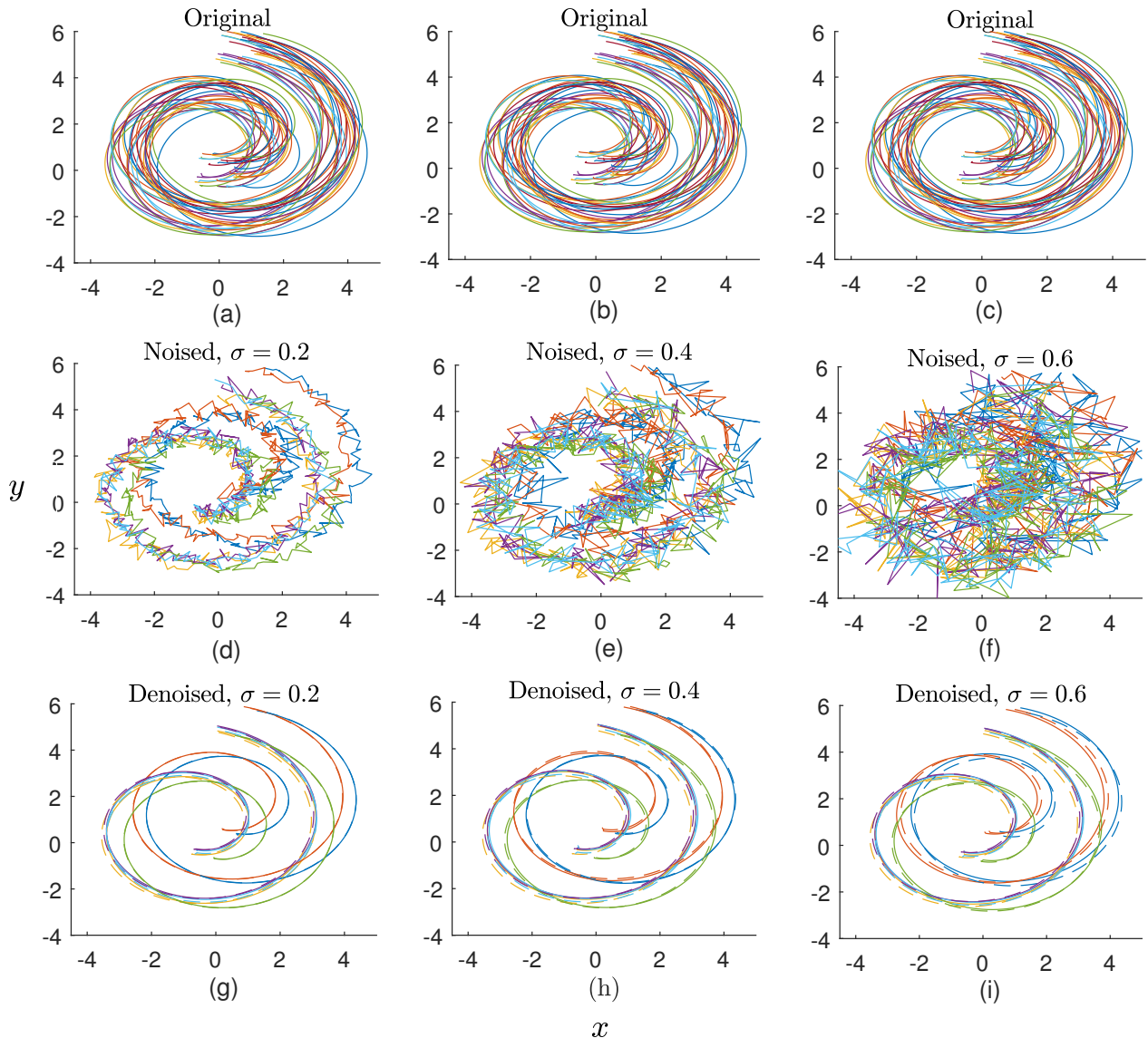}
\caption{Reconstruction of noisy collective motion trajectories using RNNs. (a-c) An original trajectory dataset of 30 agents is  generated using a generalized version of the Vicsek model (all the three datasets in a-c are identical). We make three copies of the original dataset by corrupting it with three levels of Gaussian noises with $\sigma$ = .2, .4, and .6 of $\mathcal{N}(0, \sigma^2)$ where (d-f) show 20\% of the noisy trajectories chosen for the testing. 80\% of the trajectories of each dataset are used for the training of the RNNs, while the rest is used for the testing. We set the parameters of the RNNs as follows: 201 for the sequential length, two for both the input size and the output size, 0.0005 for the learning rate, two-stacked RNN with 64 number of hidden nodes in each, 6548-14032 for the number of epochs based on the noise level, tanh for the activation function. We train each RNN with the noisy trajectories as the inputs and the original trajectories as the outputs. Then, we pass the noisy test trajectories in each (d-f) through the corresponding RNN and denoise them as shown in (g)-(i) where the denoised trajectories are shown in continuous lines while the corresponding original trajectories are shown in dashed lines with the same colors. Root mean square errors of the reconstructions are $0.05, 0.08, 0.10$ for the noise levels $\sigma$ = .2, .4, and .6 of $\mathcal{N}(0, \sigma^2)$, respectively.}
\label{fig:noiseFixAgent}
\end{figure*}

We generate 30 agents using Eqns.~\eqref{eqn:modVicsek1}, \eqref{eqn:modVicsek}, \eqref{eqn:modVicsek2}, and \eqref{eqn:rotmat2} with 201 for the time-steps ($T$), a rectangular domain with periodic boundary conditions. We set two for the radius of interaction ($r_d$), 0.05 for the speed of the particles ($\bs{v}^{(t)}_i$ for all $t$ and $i$), 0.05 for the noise on the orientation ($\epsilon^{(t)}_i$ for all $t$ and $i$), and one for the time-step size ($\delta$), see the first row of Fig.~\ref{fig:noiseFixAgent}. Then, we impose three noise levels sampled from the Gaussian distribution $\mathcal{N}(0, \sigma^2)$, where $\sigma$ = .2, .4, and .6, into the original dataset and make three copies of that, see the second row of Fig.~\ref{fig:noiseFixAgent}.

Similar to the Lorenz experiment, in order to avoid both the over-fitting and the semi-convergence of RNNs, the best number of training epochs was chosen by training an RNN with a fixed, but big, number of epochs to find the epoch number that gives the minimum loss, and then training another RNN with that many epochs. Based on the noise level of the collective motion trajectories, the best number of epochs is between 6548 and 14032. We tested both the hyperbolic tangent and Relu as activation functions and found that the loss with hyperbolic tangent is smaller than the one with ReLu. Thus, we set the activation function of the RNN to a hyperbolic tangent. Similar to Lorenz experiment, to find the optimal parameter values for the learning rate, the number of RNNs, the number of hidden nodes, first, we discretize the parameter space, then, we train and test the RNN with those discretized parameter combinations, finally, we select the best values for the parameters that give the least reconstruction error. Thus, the best learning rate is 0.0005, the best number of RNNs is two, and the best number of nodes is 64. Moreover, the other parameters of the RNN are set as follows: a batch to one trajectory, sequential length to 201 since each trajectory has 201 time-steps, input size to two since each trajectory is two-dimensional, output size to two since the output is also two-dimensional.

For each noise level, we train one RNN with 80\%, i.e., 24, of the trajectories, and test with the rest of them. We input each trajectory $\bs{x}_n$ into the RNN one at a time and generate the corresponding output of RNN, denoted as $\hat{\bs{y}}_n$. We compute the reconstruction error between $\hat{\bs{y}}_n$ and the label $\bs{y}_n$ using Eqn.~\eqref{eq:rmse}. We use BTT to compute the derivatives of Eqn.~\eqref{eq:mse} with respect to the weights and update the weights using ADAM as explained in Sec.~\ref{sec:opt}. The remaining 20\%,  i.e., six, of the noisy trajectories are fed into the corresponding trained RNN one trajectory at a time and obtain the reconstructions, $\hat{\bs{y}}_n$, see the third row of Fig.~\ref{fig:noiseFixAgent}. Root mean square errors of the reconstructions, computed using Eqn.~\eqref{eq:rmse}, are $0.05, 0.08, 0.10$ for the noise levels $\sigma$ = .2, .4, and .6 of $\mathcal{N}(0, \sigma^2)$, respectively. This test justifies that the RNN's performance decreases as the noise level increases. 

\subsection{Spiky time series of rainfall-runoff}
We use RNNs to forecast on hydrology data, especially rainfall-runoff time series which tends to be spiky. Then, compare the RNN's results with another forecast produced by a widely used hydrological model named as GR4J available in Ref.~\onlinecite{vaze2011guidelines}. Hydrology data, especially, rainfall-runoff, is highly volatile with frequent spikes so that the hydrologists have been struggling over the past decades to improve modeling of them \cite{jayathilake2021assessing}. GR4J is a daily lumped rainfall-runoff model proposed to understand hydrologic catchments' behaviors \cite{vaze2011guidelines}. This model is characterized by two independent variables, namely, precipitation (denoted by $p$) and potential evapotranspiration  (denoted by PET),  dependent variable, namely, streamflow (denoted by $q$), and four parameters, namely, the capacity of the production store (denoted by $\nu_1$), ground exchange coefficient (denoted by $\nu_2$), the capacity of the nonlinear routing store (denoted by $\nu_3$), and unit hydrograph time base (denoted by $\nu_4$). The GR4J model was used in combination with the degree-day snowmelt module available in Ref.~\onlinecite{kollat2012multiobjective}, to account for snow in northern latitudes. This degree-day snowmelt module has four parameters, namely, threshold temperature (denoted by TT), degree-day factor (denoted by CFMAX), refreezing factor (denoted by CFR), and water holding capacity of the snowpack (denoted by CWH). GR4J requires tuned values for eight parameters in order to produce accurate forecasts. We use 10 years of recent precipitation, evapotranspiration, and streamflow data from three arbitrarily chosen U.S. catchments, that we denote as Site 1 (MOPEX ID is 2365500), Site 2 (MOPEX ID is 2492000), and Site 3 (MOPEX ID is 7378000), of the model parameter estimation Experiment (MOPEX) database \cite{duan2006model}. 

\begin{figure*}[htp]
\centering
\includegraphics[width=.9\linewidth]{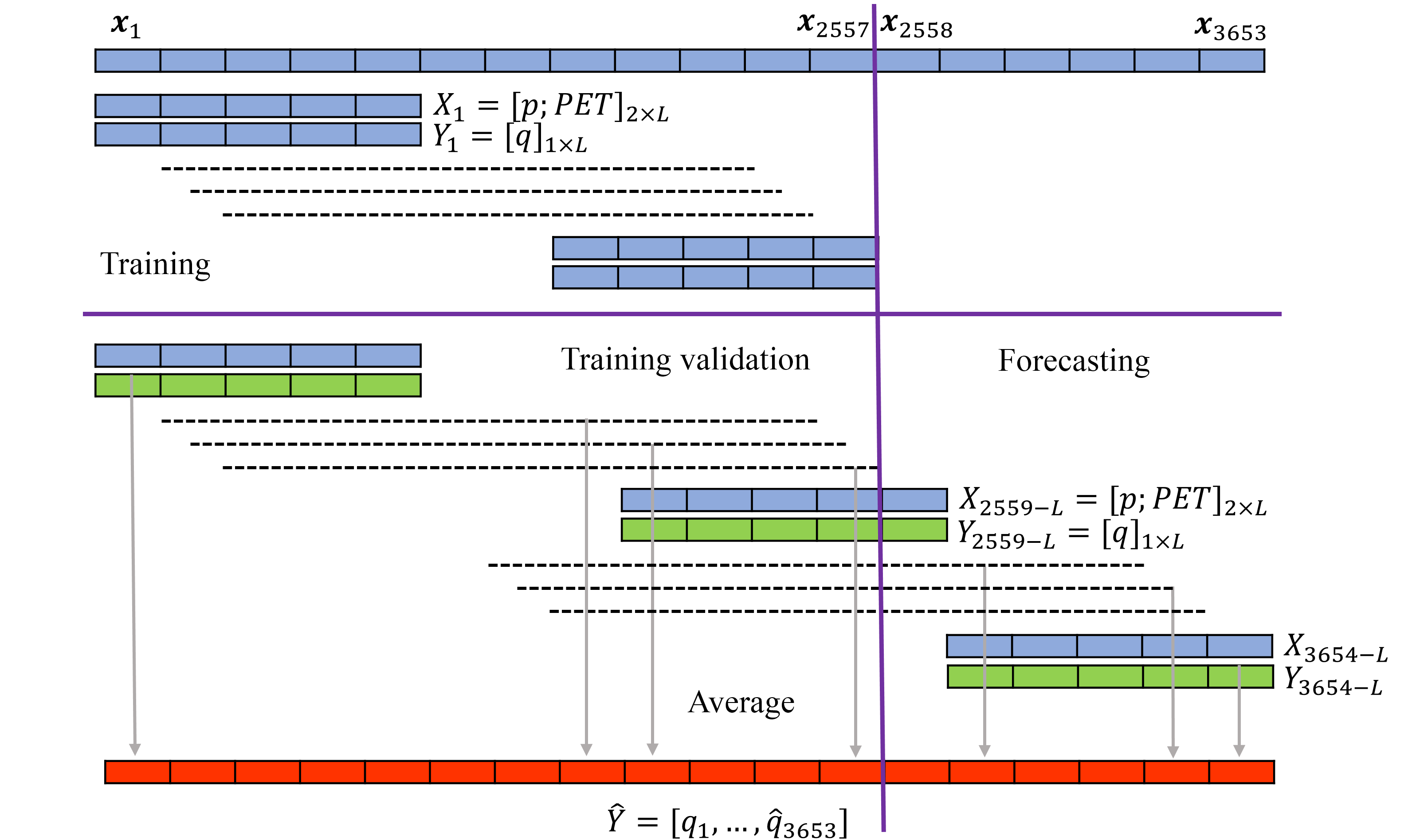}
\caption{Training and forecasting RNNs with recent 10 years of MOPEX data using a sliding window approach. We train an RNN with sliding windows of precipitation ($p$), potential evapotranspiration (PET), and streamflow  ($q$) where the first seven years, i.e., 2557 days, are used for the training to make forecasting for the last three years, i.e., 1096. We make overlapping windows of the data such that each window is $L$ time-steps long for some $L<2557$ where this figure shows a special case of $L = 5$. Here, each $\bs{x}_t$ for $t\in[1,3653]$, is a vector of $p$, PET, and $q$ so that $p$ and PET are used as the inputs to the RNN and $q$ is used as the labels during the training process, see the blue color vector pairs. Specifically, for the first iteration, we train the RNN with $X_1=[p; PET]_{2\times T}$ as the input and $Y_1=[q]_{1\times T}$ as the labels.
 Similarly, for all $t\in [2,2558-L]$, we train the same RNN with all such $X_t$'s and $Y_t$'s.  For $t\in[1,3654-L]$, we pass each $X_t$ through the RNN and obtain the modeled streamflow windows $Y_t$, see the blue-green color pairs. For each $t\in[1,3653]$, the corresponding streamflow is computed as the average of the streamflows between overlapped windows. The modeled streamflow for $t\in[1,2557]$ is used to validate the training performance of the RNN and the modeled streamflow for $t\in[2558, 3653]$ is considered as the streamflow forecast, see the red color.}
\label{fig:averageWin}
\end{figure*}

From the 10 recent years, we use the first seven years, i.e., 2557 days, for the training and the last three years, i.e., 1096 days, for the forecasting. For all three hydrological sites, we follow the same training and forecasting routings and use the same parameter values in RNNs. Similar to the previous experiments, we utilize hyperbolic tangent as the activation function and 0.001 as the learning rate. We use three-stacked RNNs with 512 nodes at each hidden layer. We utilize a special training and forecasting approach here using sliding windows as it aids sampling of adequate amount of training data just from the given time series. We make overlapping windows of the data such that each window is $L$ days long for some $L<2557$ where two consecutive windows are misaligned by only one day, see Fig.~\ref{fig:averageWin} for a special case of $L = 5$. Since we vary window size in this experiment such that $L=[5, 15, 30, 45, 60, 10, 240, 365]$, we train one RNN for each case of $L$. Thus, the sequential length parameter of the RNNs is set to the corresponding value for $L$ one at a time. Since RNN's input is two one-dimensional input vectors and the output is a single one-dimensional vector, the input size parameter is set to two and the output size parameter is set to one. Similar to the previous experiments, in order to avoid both the over-fitting and the semi-convergence of the RNN, the best number of training epochs was chosen by running the RNN with a fixed, but big, number of epochs to find the epoch number that gives the minimum loss, and then re-running the RNN with that many epochs to get the results. Based on $L$, the best number of epochs is between 99178 and 99994. We set 0.001 for the learning rate, 99200 for the number of epochs, tanh for the activation function. 

We denote $\bs{x}_t=[p,PET,q]$ for $t\in[1,3653]$ where $p$ and PET are used as the inputs to the RNN and $q$ is used as the labels during the training process. Specifically, for the first iteration, we train the RNN with $X_1=[p; PET]_{2\times T}$ as the input and $Y_1=[q]_{1\times T}$ as the labels. Similarly, for all $t\in [2,2558-L]$, we train the same RNN with all such $X_t$'s and $Y_t$'s,  see the blue color pairs in Fig.~\ref{fig:averageWin}.  For $t\in[1,3654-L]$, we pass each $X_t$ through the RNN and obtain the modeled streamflow windows $Y_t$, see the blue-green color pairs in Fig.~\ref{fig:averageWin}. For each $t\in[1,3653]$, the corresponding streamflow is computed as the average of the streamflows between overlapped windows. The modeled streamflow for $t\in[1,2557]$ is used to validate the training performance of the RNN and the modeled streamflow for $t\in[2558, 3653]$ is considered as the streamflow forecast, see the red color vectors in Fig.~\ref{fig:averageWin}. 

\begin{figure}[htp]
\centering
\includegraphics[width=1\linewidth]{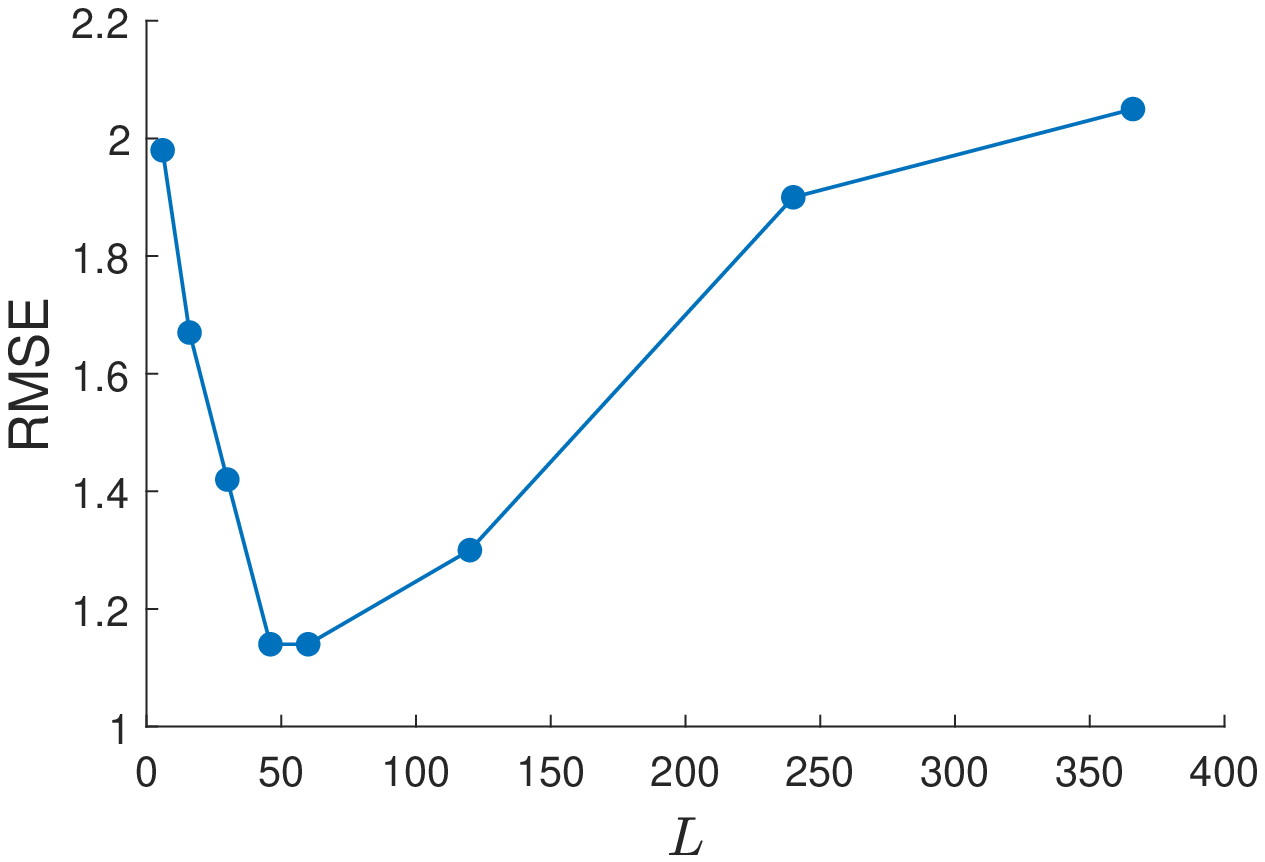}
\caption{RMSE of Site 1 (MOPEX ID is 2365500) for varying window sizes ($L$) where the RMSE is computed between the true streamflow and RNN's forecasted streamflow. We observe that the least RMSE is attained around $L=50$.}
\label{fig:winSize}
\end{figure}

\begin{figure*}[htp]
\centering
\includegraphics[width = 1\linewidth]{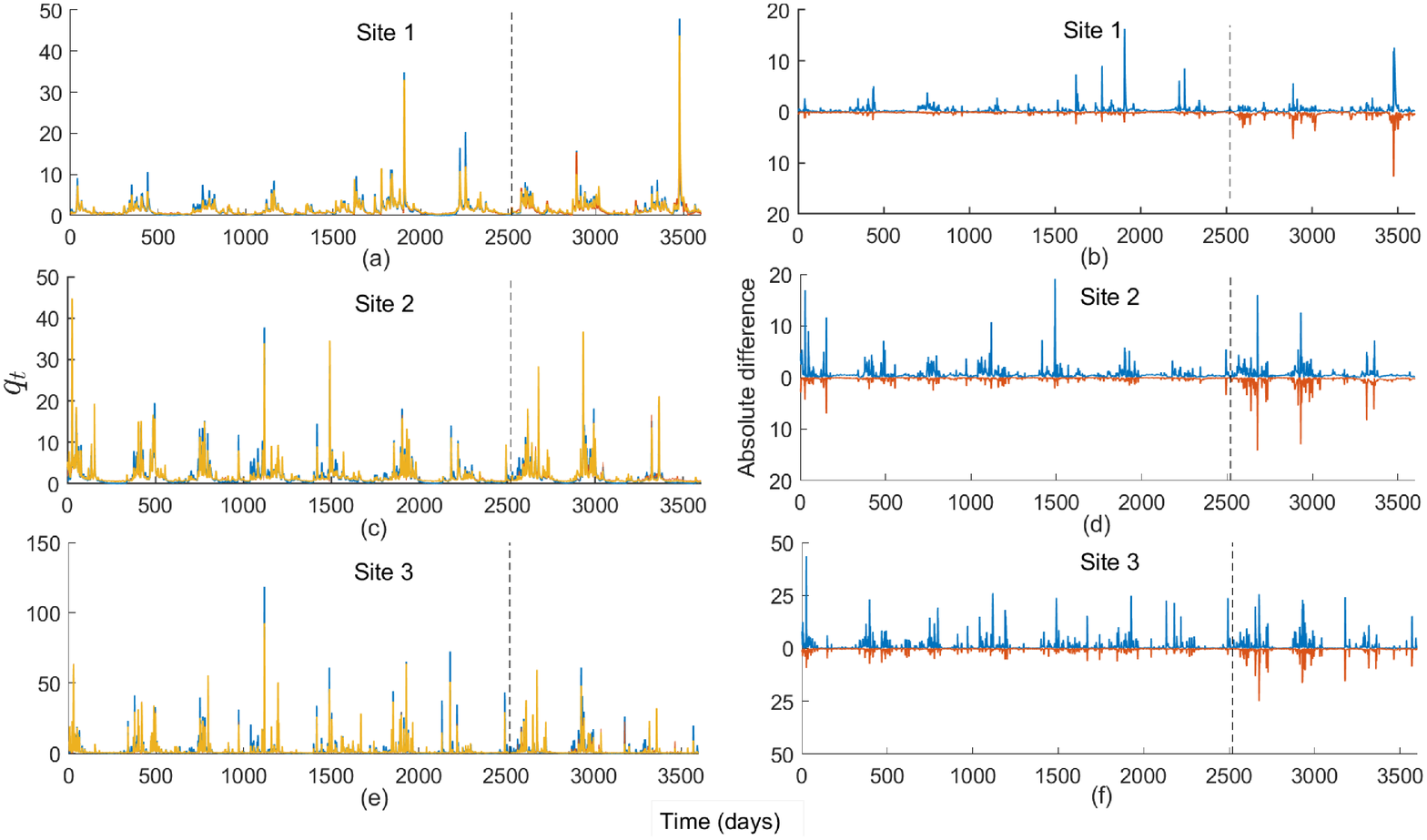}
\caption{Comparison of streamflow, denoted by $q=[q_t\vert \ \forall t]$, forecasting performance between RNN and GR4J for three hydrological sites, named as Site 1 (MOPEX ID is 2365500), Site 2 (MOPEX ID is 2492000), and Site 3 (MOPEX ID is 7378000). We use 10 recent years of data for the study where the first seven years, i.e., 2557 days, are used for training and the last three years,  i.e., 1096 days, are used for the testing. We set the same value for the same parameter between each of the RNN implemented for each site as follows: 45 for the sequential length, 2 for the input size, 1 for the output size, 0.001 for the learning rate, three-stacked RNNs, 512 for the number of hidden nodes, 99200 for the number of epochs, tanh for the activation function. We train the RNN with the first seven years of data where the inputs are precipitation, denoted by $p$, and potential evapotranspiration, denoted by PET, and the output is streamflow. Then, we pass [$p$; PET] through the trained RNN to obtain the modeled $q$ where $q$ for $t\in[1,2557]$ is used for training validation and $q$ for $t\in[2558,3653]$ is considered as the forecasting. Orange color plots in (a), (c), and (e) show true streamflow, and the red color plots therein show the RNN's modeled streamflow where the vertical black color line separates training and testing periods. Since the red color is not visible since it is underneath the orange color, we compute the absolute difference between the true and the modeled streamflows that we show by red color plots in (b), (d), and (f). We set the parameters of GR4J, namely, capacity of the production store ($\nu_1$), ground exchange coefficient ($\nu_2$), the capacity of the nonlinear routing store ($\nu_3$), and unit hydrograph time base ($\nu_4$), threshold temperature (TT), degree-day factor (CFMAX), refreezing factor (CFR), and water holding capacity of the snowpack (CWH),  to 571, -0.03, 48, 2, -0.02, 4.41, 0.36, and 0.23, respectively, for all the three sites and calibrate it. Then, we model the streamflow similar to that of RNN for the entire period where the blue color plots in (a), (c), and (e) show the GR4J's modeled streamflows and the blue color plots in (b), (d), and (f) show the absolute differences between true and modeled streamflows. The RMSE between forecasted and observed streamflows of the Sites 1, 2, and 3 are 0.9, 1.1, and 2.0 for RNNs, respectively, and are 1.0, 1.3, and 2.7 for GR4J, respectively.}
\label{fig:hydrologyDiff}
\end{figure*}

We compute RMSEs between the true and the forecasted streamflows using Eqn.~\eqref{eq:rmse} for each widow size $L$. Fig.~\ref{fig:winSize} shows RMSEs versus window size for Site 1 where we observe the least RMSE at $L=45$. The orange color plots of Figs.~\ref{fig:hydrologyDiff}(a), (c), and (e) show the true streamflows and the red color plots of those figures show the modeled streamflow by RNNs, and the red color plots in Figs.~\ref{fig:hydrologyDiff}(b), (d), and (f) show the absolute difference between the true and the RNN's modeled streamflows. The RMSEs, computed using Eqn.~\eqref{eq:rmse}, between the true and forecasted streamflows are 0.9, 1.1, and 2.0 for Sites 1, 2, and 3, respectively. 

The eight parameters of GR4J have the recommended ranges as follows: $\nu_1\in [0,1000]$, $\nu_2\in [-5,5]$, $\nu_3\in [0,300]$, $\nu_4\in [0.5,5]$, TT $\in [-3,3]$ , CFMAX $\in [0,20]$, CFR $\in [0,1]$, and CWH $\in [0.0.8]$ \cite{kollat2012multiobjective, vaze2011guidelines, jayathilake2021assessing}. We uniformly discretize GR4J's parameter space of the eight parameters into 50000 points where each point represents a combination of parameter values. We run GR4J with all of the parameter combinations for first seven years of the data where the best parameter combination with respect to RMSE is $\nu_1$ = 571, $\nu_2$ = -0.03, $\nu_3$ = 48, $\nu_4$ = 2, TT = -0.20, CFMAX = 4.41, CFR = 0.36, and CWH = 0.23, for the Site 1. We set the same values for the parameters of Site 2 and Site 3 since their optimum parameter values are either the same or sufficiently close to those of Site 1. We calibrate the GR4J model with the data for each site and model the streamflow $q_t$ for $t\in [1,3653]$ that we show by blue color plots in Figs.~\ref{fig:hydrologyDiff} (a), (b), and (c). The absolute difference of the true and GR4J's modeled streamflows are shown by the blue color plots in Figs.~\ref{fig:hydrologyDiff} (b), (d), and (f). The RMSEs between the true and the forecasted streamflows are 1.0, 1.3, and 2.7 for the sites 1, 2, and 3, respectively. We observe that the RMSEs associated with the RNN's modeled streamflows are less than that of the Gr4J's modeled streamflows. GR4J is a widely used streamflow modeling scheme\cite{jayathilake2020understanding}; however, the performance of RNN in forecasting of streamflow is capable of exceeding the performance of GR4J.

Table~\ref{tab:parameter} summarizes the parameter values used in RNNs for each of the dataset. 
\begin{table*}[htp]
\caption{Parameter values used in the RNNs trained for each datatset.}
\begin{center}
\begin{tabular}{|p{3.5cm}||p{1.4cm}|p{1.4cm}|p{1.4cm}|p{1.4cm}|p{1.4cm}|p{1.9cm}|p{1.9cm}|p{1.4cm}|}
\hline
Dataset  &Sequantial length&Input size&Output size&Learning rate &Num. of RNNs &Num. of hidden nodes&Num. of epochs&Activation function\\ 
\hline
\hline
Lorenz system with a formulation error                                   & 5000 & 3 & 3&0.01 & 3 & 128&179-930&\text{ReLu}\\
\hline
Noised trajectories of collective motion                                            & 201 & 2 & 2 &0.0005 & 2 & 64&6548-14032&\text{tanh}\\ 
\hline
Spiky time series of rainfall-runoff      & 5--365 & 2 & 1&0.001& 4 & 512&99178-99994&\text{tanh}\\
\hline
\end{tabular}\label{tab:parameter}
\end{center}
\end{table*}

\section{Conclusion}\label{sec:concl}
In this paper, we have utilized RNNs for three diverse tasks, namely, correction of a formulation error, reconstruction of corrupted particle trajectories, and forecasting of streamflow data, in three diverse fields of dynamical systems, namely, ODEs, collective motion, and hydrological modeling, respectively. The traditional approaches to solve spatiotemporal dynamical systems may not capture non-linear and complex relationships in the data since they are mostly either linear or non-linear model-based. Such spatiotemporal nonlinear systems can effectively be learned by nonlinear and model-free machine learning models such as RNNs with both high robustness and better accuracy as RNNs possess internal memory that enhance the learning process. 

We used the ODE model Lorenz system to produce two sets of data such that one consist the solutions for a correct model whereas the other consist the solution for an erroneous Lorenz model. We observed that the trained RNN was capable of producing the corresponding correct solution when the incorrect solution was provided. This approach can be used for eliminating errors associated with ODEs solvers even when only an imprecise closed-form solution is available. Future work in this context maybe using a further advanced neural network tool such as LSTM presented in Ref.~\onlinecite{hochreiter1997long} as it is composed of three gates that regulate a better information flow through the unit. Such an approach enhances the learning process of spatiotemporal data so that it better aid in eliminating the formulation error even for longer corrupted orbits.

We used a generalized version of the classic Vicsek model to generate trajectories imitating a spiral collective motion scenario. We diagnosed the influence of noise contamination on RNN's denoising performance. Such noise contamination in collective motion trajectories can occur due to the causes like less precision in video cameras and less robustness of the tracking methods \cite{gajamannage2020patch}. The results show that RNNs are capable of denoising the noisy trajectories while ensuring the pattern of the underlying collective motion. Such an ANN-based denoising method can be integrated into multi-object tracking methods to generate trustworthy results. As future work, we are planning to use RNN's denoised trajectories for the subsequent collective motion analysis tasks such as phase transition detection as we presented in Ref.~\onlinecite{gajamannage2017detecting}. 
     
Streamflow modeling is a sophisticated task in the field of hydrology due to its high volatility and frequent spikes. GR4J is a popular tool for streamflow modeling since it has shown trustworthy streamflow modeling performance. However, we use RNNs to model streamflow with a sliding window approach. The sliding window approach is capable of generating a big amount of homogeneous data from a given time series. Moreover, the sliding window approach can be used to find the most effective sequential length of the data samples that captures the most features of the given time series. Thus, we trained an RNN with sliding windows of the empirically tested best length and then model streamflow for the entire time duration. We observed that the RNN's training process involves less reconstruction error than that of GR4J. The reason for that is RNN is model-free whereas GR4J is model-based, so RNN has more tendency to adjust to any natural system than that of Gr4J. The observation that the streamflow forecasting of RNN is better than that of GR4J justifies the best performance of RNN in real-life time series forecasting. We are planning to model streamflow using other well-known models such as HBV \cite{bergstrom1995hbv} and Simhyd \cite{chiew2002application}, and then compare those results with that of RNNs. 

RNNs are capable of learning spatiotemporal characteristics of data derived from dynamical systems that we used for three such applications representing three fields. The observations based on our broad analysis justify that an RNN is an effective machine learning tool in learning a vast array of dynamical systems. The internal memory module in RNNs grants this ability which is not available in most of the spatiotemporal dynamic analyzing tools.

\section*{Acknowledgments}
The authors would like to thank the Google Cloud Platform for granting Research Credit to access its GPU computing resources under the project number 397744870419. E.~M.~Bollt has received funding from the Army Research Office (No. N68164-EG), the National Institutes of Health (NIH-CRCNS) and also the Defense Advanced Research Projects Agency.

\section*{References}

%

\end{document}